\documentclass[12pt]{amsart}
\usepackage{mathrsfs}
\usepackage{pifont}
\usepackage{amsfonts}
\usepackage{amsthm, amsmath, amssymb}
\usepackage{longtable}
\usepackage{multirow,makecell}
\usepackage[labelfont=bf,textfont=bf]{caption}
\usepackage{tabu}
\usepackage{booktabs}
\usepackage{rotating} 
\usepackage{tabularx}
\usepackage{color}

\renewcommand{\proof}{\noindent{\it Proof.\ \ }}
\def\qed{\ifmmode\square\else\nolinebreak\hfill
$\Box$\fi\par\vskip12pt}

\def\la{\langle} \def\ra{\rangle} 
\def\div{\,\big|\,} 

    \def\ZZ{\mathbb Z}

\def\diam{{\sf diam}}

\def\K{{\sf K}}

\def\Ga{{\it \Gamma}} \def\Sig{{\it \Sigma}} 
\def\Ome{{\it \Omega}}

\def\Cay{{\sf Cay}} \def\Cos{{\sf Cos}}  
\def\Aut{{\sf Aut}}

   \def\K{{\sf K}}

\def\a{\alpha} \def\b{\beta} \def\d{\delta}\def\g{\gamma}

\def\A{{\sf A}}

\newtheorem{theorem}{Theorem}[section]%
\newtheorem{lemma}[theorem]{Lemma}%
\newtheorem{corollary}[theorem]{Corollary}%
\newtheorem{proposition}[theorem]{Proposition}%
\newtheorem{question}[theorem]{Question}%

\begin{document}

\title[Two-distance transitive normal Cayley graphs]
{Two-distance transitive normal Cayley graphs}

\thanks{2010 MR Subject Classification 05C25, 05E18, 20B25.}
\thanks{Corresponding author: Yan-Quan Feng. }
\thanks{The work was supported by the National Natural Science Foundation of China (11731002) and the
111 Project of China (B16002).}

\author[Huang, Feng, and Zhou]{ Jun-Jie Huang, Yan-Quan Feng, Jin-Xin Zhou}
\address{Jun-Jie Huang, Yan-Quan Feng, Jin-Xin Zhou\\
Department of Mathematics\\
Beijing Jiaotong University\\
Beijing \\
100044, P. R. China}
\email{20118006@bjtu.edu.cn(J.-H. Huang), yqfeng@bjtu.edu.cn (Y.-Q. Feng), \linebreak
  jxzhou@bjtu.edu.cn (J.-X. Zhou)}
\maketitle

\begin{abstract}
In this paper, we construct an infinite family of normal Cayley graphs,
which are $2$-distance-transitive but neither distance-transitive nor $2$-arc-transitive. This
answers a question raised by Chen, Jin and Li in 2019 and corrects a claim in a literature given by Pan, Huang and Liu in 2015.
\end{abstract}

\qquad {\textsc k}{\scriptsize \textsc {eywords.}} {\footnotesize
Cayley graph, $2$-distance-transitive graph, simple group.}

\section{Introduction}

In this paper, all graphs are finite, simple, and undirected.
For a graph $\Ga$, let $V(\Ga), E(\Ga), A(\Ga)$ or $\Aut(\Ga)$ denote its vertex set,
edge set, arc set and its full automorphism group, respectively.
The graph $\Ga$ is called $G$-{\it vertex-transitive}, $G$-{\it edge-transitive} or $G$-{\it arc-transitive}, with $G\le\Aut(\Ga)$, if $G$ is transitive on $V(\Ga),E(\Ga)$ or $A(\Ga)$ respectively, and $G$-{\em semisymmetric}, if $\Ga$ is $G$-edge-transitive but not $G$-vertex-transitive. It is easy to see that a $G$-semisymmetric graph $\Ga$ must be bipartite such that $G$ has two orbits, namely the two parts of $\Ga$, and the stabilizer $G_u$ for any $u\in V(\Ga)$ is transitive the neighbourhood of $u$ in $\Ga$. An $s$-{\it arc} of $\Ga$ is a sequence $v_0,v_1,\dots,v_s$ of $s+1$ vertices of $\Ga$ such that $v_{i-1},v_i$ are adjacent for $1\le i\le s$ and $v_{i-1}\ne v_{i+1}$ for $1\le i\le s-1$. If $\Ga$ has at least one $s$-arc and $G\le\Aut(\Ga)$ is transitive on the set of $s$-arcs of $\Ga$, then $\Ga$ is called {\it $(G,s)$-arc-transitive}, and $\Ga$ is said to be {\it $s$-arc-transitive} if it is $(\Aut(\Ga),s)$-arc-transitive.

For two vertices $u$ and $v$ in $V(\Ga)$, the {\it distance} $d(u,v)$ between $u$ and $v$ in $\Ga$ is the smallest length of paths between $u$ and $v$, and the {\it diameter} $\diam(\Ga)$ of $\Ga$ is the maximum distance occurring over all pairs of vertices.
For $i=1,2,\cdots,\diam(\Ga)$, denote by $\Ga_i(u)$ the set of vertices at distance $i$ with vertex $u$ in $\Ga$. A graph $\Ga$ is called {\it distance transitive} if, for any vertices $u,v,x,y$ with $d(u,v)=d(x,y)$,
there exists $g\in\Aut(\Ga)$ such that $(u,v)^g=(x,y)$.
The graph $\Ga$ is called {\it $(G,t)$-distance-transitive} with $G\leq\Aut(\Ga)$ if, for each $1\leq i\leq t$, the group $G$ is transitive on the ordered pairs of form $(u,v)$ with $d(u,v)=i$, and $\Ga$ is said to be {\em $t$-distance-transitive} if it is $(\Aut(\Ga),t)$-distance-transitive.

Distance-transitive graphs were first defined by Biggs and Smoth in \cite{BS}, and
they showed that there are only 12 trivalant distance-transitive graphs.
Later, distance-transitive graphs of valencies $3$, $4$, $5$, $6$ and $7$
were classified in~\cite{BS,FII84,GP85,GP86,GP87},
and a complete classification of distance-transitive graphs with symmetric or alternating groups of automorphisms was given by Liebeck, Praeger and Saxl~\cite{LPS87}.
The $2$-distance-transitive but not 2-arc-transitive graphs of valency at most 6 were classified in \cite{CJS17,JL14},
and the $2$-distance-primitive graphs (a vertex stabilizer of automorphism group is primitive on
both the first step and the second step neighbourhoods of the vertex) with prime valency were classified in~\cite{JHL}.
By definition, a $2$-arc-transitive graph is $2$-distance-transitive, but a $2$-distance-transitive graph may not be $2$-arc-transitive; a simple example is the complete multipartite graph $\K_{3,2}$.
Furthermore, Corr, Jin and Schneider~\cite{CJS} investigated properties of a connected $(G,2)$-distance-transitive but not $(G,2)$-arc-transitive graph of girth $4$,
and they applied the properties to classify such graphs with prime valency. For more information about $2$-distance-transitive graphs, we refer to \cite{DGLP12,DGLP13}.

For a finite group $G$ and a subset $S\subseteq G \setminus\{1\}$ with $S=S^{-1}:=\{s^{-1}\mid s\in S\}$, the {\it Cayley graph} $\Cay(G,S)$ of the group $G$ with respect to $S$ is the graph with vertex set $G$ and with two vertices $g$ and $h$ adjacent if $hg^{-1}\in S$.
For $g\in G$, let $R(g)$ be the permutation of $G$ defined by $x\mapsto xg$ for
all $x\in G$. Then $R(G):=\{R(g)\mid g\in G\}$ is a regular group of automorphisms of  $\Cay(G,S)$.
It is known that a graph $\Ga$ is a Cayley graph of $G$ if and only if
$\Ga$ has a regular group of automorphisms on the vertex set which is isomorphic to $G$;
see \cite[Lemma 16.3]{Biggs} and \cite{Sabidussi}.
A Cayley graph $\Ga=\Cay(G,S)$ is called {\it normal} if $R(G)$ is a normal subgroup of $\Aut(\Ga)$. The study of normal Cayley graphs was initiated by Xu \cite{Xu98} and has been investigated under various additional conditions; see \cite{DJLP,Praeger99}.

There are many interesting examples of arc-transitive graphs and
$2$-arc-transitive graphs constructed as normal Cayley graphs. However, the status for $2$-distance-transitive graphs is different. Recently, $2$-distance-transitive circulants
were classified in \cite{CJL}, where the following question was proposed:
\begin{question}{\rm(\cite[Question 1.2]{CJL})}\label{question}
Is there a normal Cayley graph which is 2-distance-transitive,
but neither distance-transitive nor 2-arc-transitive?
\end{question}

In this paper, we answer the above question by constructing an infinite family of such graphs.

\begin{theorem}\label{Thm-1} For an odd prime $p$, let $G=\la a,b,c \mid a^p=b^p=c^p=1,[a,b]=c,[c,a]=[c,b]=1\ra$ and $S=\{a^i,b^i\mid1\leq i\leq p-1\}$.
Then $\Cay(G,S)$ is a $2$-distance-transitive normal Cayley graph that is neither distance-transitive nor $2$-arc-transitive.
\end{theorem}

Applying this theorem, we can obtain the following corollary.

\begin{corollary} \label{cor}
Under the notation given in {\em Theorem~\ref{Thm-1}}, let $\Cos(G, \la a\ra, \la b\ra)$ be the graph with vertex set $\{\la a\ra g\ |\ g \in G\} \cup \{\la b\ra h\ |\ h \in G\}$ and with edges all these coset pairs $\{\la a\ra g, \la b\ra h\}$ having non-empty intersection in $G$. Then $\Cay(G,S)$ is the line graph of $\Cos(G, \la a\ra, \la b\ra)$, and $\Cos(G, \la a\ra, \la b\ra)$ is $3$-arc-transitive.
\end{corollary}

The graph $\Cos(G, \la a\ra, \la b\ra)$ was first constructed in \cite{PHL} as a regular cover of $\K_{p,p}$, where  it is said that $\Cos(G, \la a\ra, \la b\ra)$ is $2$-arc-transitive in \cite[Theorem~1.1]{PHL}, but not $3$-arc-transitive generally for all odd primes $p$ in a remark after \cite[Example~4.1]{PHL}. However, this is not true and Corollary~\ref{cor} implies that $\Cos(G, \la a\ra, \la b\ra)$ is always $3$-arc-transitive for each odd prime $p$. In fact, $\Cos(G, \la a\ra, \la b\ra)$ is $3$-arc-regular, that is, $\Aut(\Cos(G, \la a\ra, \la b\ra))$ is regular on the set of $3$-arcs of $\Cos(G, \la a\ra, \la b\ra)$.

\section{Preliminaries}

In this section we list some preliminary results used in this paper. The first one is the well-known orbit-stabilizer theorem (see \cite[Theorem 1.4A]{Dixon}).

\begin{proposition}\label{orbit-stab}
Let a group $G$ has a transitive action on a set $\Omega$ and let $\a\in \Omega$. Then $|G|=|\Omega||G_\a|$.
\end{proposition}

The well-known Burnside $p^aq^b$ theorem was given in~\cite[Theorem 3.3]{Gorenstein}.

\begin{proposition}\label{Burnside}
Let $p$ and $q$ be primes and let $a$ and $b$ be positive integers. Then a group of order $p^aq^b$ is soluble.
\end{proposition}

The next proposition is an important property of a non-abelian simple group acting transitively on a set with cardinality a prime-power, and we refer to \cite[Corollary 2]{Guralnick} or \cite[Proposition 2.4]{WFZ}.

\begin{proposition}\label{2-trans}
Let $T$ be a nonabelian simple group acting transitively on a set $\Ome$ with cardinality a $p$-power for a prime $p$. If $p$ does not divide the order of a point-stabilizer of $T$, then $T$ acts $2$-transitively on $\Ome$.
\end{proposition}

Let $\Ga=\Cay(G,S)$ be a Cayley graph of a group $G$ with respect to $S$. Then $R(G)$ is a regular subgroup of $\Aut(\Ga)$, and $\Aut(G,S):=\{\a\in\Aut(G)\mid S^\a=S\}$ is also a subgroup of $\Aut(\Ga)$, which fixes $1$. Furthermore, $R(G)$ is normalized by $\Aut(G,S)$, and hence we have a semiproduct $R(G)\rtimes\Aut(G,S)$, where $R(g)^\a=R(g^\a)$ for any $g\in G$ and $\a\in \Aut(G,S)$.
Godsil~\cite{Godsil} proved that the semiproduct $R(G)\rtimes\Aut(G,S)$ is in fact the normalizer of $R(G)$ in $\Aut(\Ga)$. By Xu~\cite{Xu98}, we have the following proposition.

\begin{proposition}\label{Aut-Cay}
Let $\Ga=\Cay(G,S)$ be a Cayley graph of a finite group $G$ with respect to $S$, and let $\A=\Aut(\Ga)$.
Then the following hold:
\begin{itemize}
  \item [(1)] $N_\A(R(G))=R(G)\rtimes\Aut(G,S)$;
  \item [(2)] $\Ga$ is a normal Cayley graph if and only if $\A_1=\Aut(G,S)$, where $\A_1$ is the stabilizer of $1$ in $\A$.
\end{itemize}
\end{proposition}

Let $\Ga$ be a $G$-vertex-transitive graph, and let $N$ be a normal subgroup of $G$.
The {\it normal quotient graph} $\Ga_N$ of $\Ga$ induced by $N$ is defined to be the graph with vertex set the orbits of $N$ and with two orbits $B,C$ adjacent if some vertex in $B$ is adjacent to some vertex in $C$ in $\Ga$. Furthermore, $\Ga$ is called a {\it normal $N$-cover} of $\Ga_N$ if $\Ga$ and $\Ga_N$ have the same valency.

\begin{proposition}\label{Praeger}
Let $\Ga$ be a connected $G$-vertex-transitive graph and let $N$ be a normal subgroup of $G$. If $\Ga$ is a normal $N$-cover of $\Ga_N$, then the following statements hold:
\begin{itemize}
\item[(1)] $N$ is semiregular on $V\Ga$ and is the kernel of $G$ acting $V(\Ga_N)$, so $G/N\le\Aut(\Ga_N)$;
\item[(2)] $\Ga$ is $(G,s)$-arc-transitive if and only if $\Ga_N$ is $(G/N,s)$-arc-transitive;
\item[(3)] $G_{\a}\cong(G/N)_{\d}$ for any $\a\in V\Ga$ and $\d\in V(\Ga_N)$.
\end{itemize}

In particular, the above results hold if we replace the assumption that $\Ga$ is a normal $N$-cover of $\Ga_N$ by the following assumption: $\Ga$ is $G$-arc-transitive with a prime valency and $N$ has at least three orbits.
\end{proposition}

Proposition~\ref{Praeger} was given in many papers by replacing the condition that $\Ga$ is a normal $N$-cover of $\Ga_N$ by one of the following assumptions:
(1) $N$ has at least $3$-orbits and $G$ is $2$-arc-transitive (see \cite[Theorem 4.1]{Praeger92});
(2) $N$ has at least $3$-orbits, $G$ is arc-transitive and $\Ga$ has a prime valency (see \cite[Theorem 2.5]{PY});
(3) $N$ has at least $3$-orbits and $G$ is locally primitive (see \cite[Lemma 2.5]{Li-Pan}).
The first step for these proofs is to show that for any two vertices $B,C\in V(\Ga_N)$,
the induced subgraph $[B]$ of $B$ in $\Ga$ has no edge and if $B$ and $C$ are adjacent in $\Ga_N$ then the induced subgraph $[B\cup C]$ in $\Ga$ is a matching,
which is equivalent to that $\Ga$ is a normal $N$-cover of $\Ga_N$.
Then Proposition~\ref{Praeger}~(1)-(3) follows from these proofs.



\section{Proof Theorem~\ref{Thm-1}}

For a positive integer $n$ and a prime $p$, we use $\ZZ_n$ and $\ZZ_p^r$ to denote the cyclic group of order $n$ and the elementary abelian group of order $p^r$, respectively. In this section, we always assume that $p$ is an odd prime, and denote by $\ZZ_p^*$ the multiplicative group of $\ZZ_p$ consisting of all non-zero numbers in $\ZZ_p$. Note that $\ZZ_p^*\cong\ZZ_{p-1}$. Furthermore, we  also set the following assumptions in this section:
$$G=\la a,b,c\mid a^p=b^p=c^p=1,[a,b]=c,[c,a]=[c,b]=1\ra,\ S=\{a^i,b^i\mid1\leq i\leq p-1\},$$
$$\Ga=\Cay(G,S), \ \ \A=\Aut(\Ga), \ \ N=N_\A(R(G))=R(G)\rtimes\Aut(G,S), \ \mbox{ and } \ZZ_p^*=\langle t\rangle.$$

By Proposition~\ref{Aut-Cay}, $N_\A(R(G))=R(G)\rtimes\Aut(G,S)$, and $R(g)^\delta=R(g^\delta)$ for any $R(g)\in R(G)$ and $\delta\in \Aut(G,S)$. Since $G=\la S\ra$, $\Ga$ is a connected Cayley graph of valency $2(p-1)$. Let
\begin{eqnarray*}
\begin{matrix}
\a: a \longmapsto a^t,   &  b \longmapsto b,  & c \longmapsto c^t;  \\
\b: a \longmapsto a,    & b \longmapsto b^t, &c \longmapsto c^t;\\
\g: a \longmapsto b, &  b \longmapsto a,~ &c \longmapsto c^{-1}.
\end{matrix}
\end{eqnarray*}

Since $a^t, b, c^t$ satisfy the same relations as $a,b,c$ in $G$ and $G=\langle a^t, b, c^t\rangle$, $\a$ induces an automorphism of $G$, and we still denote by $\a$ this automorphism. Similarly, $\b$ and $\g$ are also automorphisms of $G$.

\begin{lemma}\label{arc-trans}
$\Aut(G,S)=\la\a,\b,\g\ra\cong(\ZZ_{p-1}\times\ZZ_{p-1})\rtimes\ZZ_2$, and $\Ga$ is $N$-arc-transitive. Furthermore, $N$ has no normal subgroup of order $p^2$.
\end{lemma}

\proof Since $\ZZ_p^*=\langle t\rangle$, it is easy to check that $\a^{p-1}=\b^{p-1}=\g^2=1$, $\a\b=\b\a$ and $\a^\g=\b$. Thus $\la\a,\b,\g\ra\cong(\ZZ_{p-1}\times\ZZ_{p-1})\rtimes\ZZ_2$. Clearly, $\a,\b,\g\in \Aut(G,S)$.
To prove $\Aut(G,S)=\la\a,\b,\g\ra\cong(\ZZ_{p-1}\times\ZZ_{p-1})\rtimes\ZZ_2$, it suffices to show that $|\Aut(G,S)|\leq 2(p-1)^2$.

Clearly, $\langle \a,\b,\g\rangle$ is transitive on $S$, and hence $\Ga$ is $N$-arc-transitive. Since $G=\langle S\rangle$, $\Aut(G,S)$ is faithful on $S$. By Proposition~\ref{orbit-stab}, $|\Aut(G,S)|=|S||\Aut(G,S)_a|$, where $\Aut(G,S)_a$ is the stabilizer of $a$ in $\Aut(G,S)$. Note that $\Aut(G,S)_a$ fixes $a^i$ for each $1\leq i\leq p-1$. Again by Proposition~\ref{orbit-stab}, $|\Aut(G,S)_a|\leq (p-1)|\Aut(G,S)_{a,b}|$, where $\Aut(G,S)_{a,b}$ is the subgroup of $\Aut(G,S)$ fixing $a$ and $b$. Since $G=\langle a,b\rangle$, we obtain $\Aut(G,S)_{a,b}=1$, and then $|\Aut(G,S)|\leq 2(p-1)^2$, as required.

Let $H\leq N$ be a subgroup of order $p^2$. Since $R(G)$ is the unique normal Sylow $p$-subgroup of $N=R(G)\rtimes\Aut(G,S)$, we have $H\leq R(G)$, and since $|R(G):H|=p$, we have $H\unlhd
R(G)$. Note that the center $C:=Z(R(G))=\langle R(c)\rangle$ and $C\cap H\not=1$. Thus,
$C\cap H=C$ as $|C|=p$, implying $C\leq H$. Since $H/C$ is a subgroup of order $p$, and  $R(G)/C=\langle R(a)C\rangle\times\langle R(b)C\rangle\cong\ZZ_p^2$, we have $H/C=\langle R(b)C\rangle$ or $\langle R(a)R(b)^iC\rangle$ for some $0\leq i\leq p-1$. It follows that $H=\langle R(b)\rangle\times C$ or $\langle R(ab^i)\rangle\times C$ for some $0\leq i\leq p-1$.

Suppose $H\unlhd N$. Since $C$ is characteristic in $R(G)$ and $R(G)\unlhd N$, we have $C\unlhd N$.
Recall that $R(a)^\g=R(a^\g)=R(b)$. Then $(\langle R(a)\rangle\times C)^\g=\langle R(b)\rangle\times C$. This implies that both $\langle R(a)\rangle\times C$ and $\langle R(b)\rangle\times C$ are not normal in $N$. Thus, $H=\langle R(ab^i)\rangle\times C$ for some $1\leq i\leq p-1$.
Since $H\unlhd N$, we have $H^\b=H$, that is, $\langle R(ab^{ti})\rangle \times C=H^\b=H=\langle R(ab^i)\rangle\times C$. It follows that $\langle R(ab^{ti})\rangle=\langle R(ab^i)\rangle$ and then  $R(ab^{ti})=R(ab^i)$, which further implies $b^{ti}=b^i$. This gives rise to $p\div i(t-1)$, and since $(i,p)=1$, we have $t=1$, contradicting that $\ZZ_p^*=\langle t\rangle\cong \ZZ_{p-1}$. Thus, $N$ has no normal subgroup of order $p^2$.  \qed

A {\it clique} of a graph $\Ga$ is a maximal complete subgraph, and the {\it clique graph} $\Sig$ of $\Ga$ is defined to have the set of all cliques of $\Ga$ as its vertex set with two cliques adjacent in $\Sig$ if the two cliques have at least one common vertex. For a positive integer $n$, $n_p$ denotes the largest $p$-power diving $n$.

\begin{lemma}\label{clique}
The clique graph $\Sig$ of $\Ga$ is a connected $p$-valent bipartite graph of order $2p^2$, $\A$ has a faithful natural action on $\Sig$, and $\Sig$ is $R(G)$-semisymmetric and $N$-arc-transitive. Furthermore, $|\A|_p=p^3$.
\end{lemma}

\proof Recall that $G=\la a,b,c\mid a^p=b^p=c^p=1,[a,b]=c,[c,a]=[c,b]=1\ra$ and $S=\{a^i,b^i\mid1\leq i\leq p-1\}$. Then $\Ga=\Cay(G,S)$ has exactly two cliques passing through $1$, that is, the induced subgraphs of $\langle a\rangle$ and $\langle b\rangle$ in $\Ga$. Since $R(G)\leq \Aut(\Ga)$ is transitive on vertex set, each clique of $\Ga$ is an induced subgraph of the coset $\langle a\rangle x$ or $\langle b\rangle x$ for some $x\in G$. Thus, we may view the vertex set of $\Sig$ as $\{\langle a\rangle x, \langle b\rangle x \ | \ x\in G\}$ with two cosets  adjacent in $\Sig$ if they have non-empty intersection. It is easy to see that $\langle a\rangle x\cap \langle b\rangle y\not=\emptyset$ if and only if $|\langle a\rangle x\cap \langle b\rangle y|=1$, and any two distinct cosets, either in $\{\langle a\rangle x\ |\ x\in G\}$ or in $\{\langle b\rangle x\ |\ x\in G\}$, have empty intersection. Furthermore, $\langle a\rangle$ has non-empty intersection with exactly $p$ cosets, that is, $\langle b\rangle a^i$ for $0\leq i\leq p-1$. Thus, $\Sig$ is a $p$-valent bipartite graph of order $2p^2$. The connectedness of $\Sig$ follows from that of $\Ga$.

Clearly, $\A$ has a natural action on $\Sig$. Let $K$ be the kernel of $\A$ on $\Sig$. Then $K$ fixes each coset of $\langle a\rangle x$ and $\langle b\rangle x$ for all $x\in G$. Since $\langle a\rangle x\cap\langle b\rangle x=\{x\}$, $K$ fixes $x$ and hence $K=1$. Thus, $\A$ is faithful on $\Sig$ and we may let $\A\leq \Aut(\Sig)$.

Note that $R(G)$ is not transitive on $\{\langle a\rangle x, \langle b\rangle x \ | \ x\in G\}$, but transitive on $\{\langle a\rangle x\ | \ x\in G\}$ and $\{\langle b\rangle x \ | \ x\in G\}$. Furthermore, $R(\langle a\rangle)$ fixes $\langle a\rangle$ and is transitive on $\{\langle b\rangle a^i\ |\ 0\leq i\leq p-1\}$, the neighbourhood of $\langle a\rangle$ in $\Sig$, and similarly, $R(\langle b\rangle)$ fixes $\langle b\rangle$ and is transitive on the neighbourhood $\{\langle a\rangle b^i\ |\ 0\leq i\leq p-1\}$ of $\langle b\rangle$ in $\Sig$. It follows that $\Sig$ is $R(G)$-semisymmetric. Recall that $N=R(G)\rtimes\Aut(G,S)$ and $\Aut(G,S)=\langle \a,\b,\g\rangle$.
Since $a^\g=b$ and $b^\g=a$, $\g$ interchanges $\{\langle a\rangle x\ | \ x\in G\}$ and $\{\langle b\rangle x \ | \ x\in G\}$. This yields that $\Sig$ is $R(G)\rtimes\langle \g\rangle$-arc-transitive and hence $N$-arc-transitive.

Since $\Sig$ is a connected graph with prime valency $p$, we have $p^2\nmid |\Aut(\Sig)_u|$ for any $u\in V(\Sig)$, and in particular, $p^2\nmid |\A_u|$. Note that $p\div |\A_u|$. By Proposition~\ref{orbit-stab}, $|\A|=|\Sig||\A_u|=2p^2|\A_u|$. This implies that $|\A|_p=p^3$. \qed

\begin{lemma}\label{normal}
$\A=\Aut(\Ga)=R(G)\rtimes\Aut(G,S)$.
\end{lemma}

\proof By Lemma~\ref{clique}, $|\A|_p=p^3$, and since $|V(\Ga)|=p^3$ and $\A$ is vertex-transitive on $V(\Ga)$, the vertex stabilizer $\A_1$ is a $p'$-group, that is, $p\nmid |\A_1|$. To prove the lemma, by Proposition~\ref{Aut-Cay} we only need to show that $R(G)\unlhd \A$, and since $R(G)$ is a Sylow $p$-subgroup of $\A$, it suffices to show that $\A$ has a normal Sylow $p$-subgroup.

Let $M$ be a minimal normal subgroup of $\A$. Then $M=T_1\times T_2\cdots\times T_d$, where $T_i\cong T$ for each $1\leq i\leq d$ with a simple group $T$. Since $|V(\Ga)|=p^3$, each orbit of $M$ has length a $p$-power and hence each orbit of $T_i$ has length a $p$-power. It follows that $p\div |T|$. Assume that $|T|_p=p^\ell$. Then $|M|_p=p^{d\ell}$ and $d\ell=1,2$ or $3$ as $|\A|_p=p^3$.

We process the proof by considering the two cases: $M$ is insoluble or soluble.

\vskip 0.2cm
\noindent{\bf Case 1:} $M$ is insoluble.

In this case, $T$ is a non-abelian simple group. We prove that this case cannot happen by deriving contradictions. Recall that $d\ell=1,2$ or $3$.

Assume that $d\ell=1$. Then $|M|_p=p$. By Lemma~\ref{clique}, $M\unlhd \A\leq \Aut(\Sig)$, and since $|V(\Sig)|=2p^2$, $M$ has at least three orbits. Since $\Sig$ has valency $p$, Proposition~\ref{Praeger} implies that $M$ is semiregular on $V(\Sig)$ and hence $|M|\div 2p^2$. By Proposition~\ref{Burnside}, $M$ is soluble, a contradiction.

Assume that $d\ell=2$. Since $R(G)$ is a Sylow $p$-subgroup of $\A$ and $M\unlhd \A$, $R(G)\cap M$ is a Sylow $p$-subgroup of $M$ and hence $|R(G)\cap M|=|M|_p=p^2$. Since $R(G)\unlhd N$ and $M\unlhd \A$, $M\cap R(G)$ is a normal subgroup of order $p^2$ in $N$, contradicting to Lemma~\ref{arc-trans}.

Assume that $d\ell=3$. Then $(d,\ell)=(1,3)$ or $(3,1)$. Since $|M|_p=p^3=|\A|_p$, we deduce $R(G)\leq M$ and hence $M$ is transitive on $\Ga$.

For $(d,\ell)=(1,3)$, $M$ is a non-abelian simple group. Since $M_1\leq \A_1$ is a $p'$-group,  Proposition~\ref{2-trans} implies that $M$ is $2$-transitive on $\Ga$, forcing that $\Ga$ is the complete graph of order $p^3$, a contradiction.

For $(d,\ell)=(3,1)$, we have $M=T_1\times T_2\times T_3$. Then $|M|_p=p^3$, and since $M\unlhd \A$, we derive $R(G)\leq M$. By Lemma~\ref{clique} $M\leq \Aut(\Sig)$, and $\Sig$ is $R(G)$-semisymmetric. Since $M$ has no subgroup of index $2$, $M$ fixes the two parts of $\Sig$ setwise, and hence $\Sig$ is $M$-semisymmetric. Noting that $\g$ interchanges the two parts of $\Sig$, we have that $\Sig$ is  $M\langle\g\rangle$-arc-transitive. Since $\g$ is an involution, under conjugacy it fixes $T_i$ for some $1\leq i\leq 3$, say $T_1$. Then $T_1\unlhd  \langle M,\g\rangle$ and by Proposition~\ref{Praeger}, $T_1$ is semiregular on $\Sig$. This gives rise to $|T_1|\div 2p^2$, contrary to the simplicity of $T_1$.

\vskip 0.2cm
\noindent{\bf Case 2:} $M$ is soluble.

Since $p\div |M|$, we have $M=\ZZ_p^d$ with $1\leq d\leq 3$. If $d=3$ then $\A$ has a normal Sylow $p$-subgroup, as required. If $d=2$ then $M\leq R(G)\leq N$ and $N$ has a normal subgroup of order $p^2$, contrary to Lemma~\ref{arc-trans}. Thus, we may let $d=1$, and since $M\leq R(G)$ and $R(G)$ has a unique normal subgroup of order $p$ that is the center of $R(G)$, we derive that $M=\langle R(c)\rangle$.

Now it is easy to see that the quotient graph $\Ga_M=\Cay(G/M,S/M)$ with $S/M=\{a^iM,b^iM\ |\ 1\leq i\leq p-1\}$. Note that $G/M=\langle aM\rangle\times\langle bM\rangle\cong \ZZ_p^2$. Then $\Ga_M$ is a connected Cayley graph of order $p^2$ with valency $2(p-1)$, so $\Ga$ is a normal $M$-cover of $\Ga_M$. By Proposition~\ref{Praeger}, we may let $\A/M\leq \Aut(\Ga_M)$ and $\Ga_M$ is $\A/M$-arc-transitive.

Let $H/M$ be a minimal normal subgroup of $\A/M$. Then $H\unlhd \A$ and $H/M=L_1/M\times \cdots \times L_r/M$, where $L_i\unlhd H$ and $L_i/M$ ($1\leq i\leq r$) are isomorphic simple groups. Since $|\Ga_M|=p^2$, we infer $p\div |H/M|$ and similarly, $p\div |L_i/M|$. Let $|L_i/M|_p=p^s$. Then $|H/M|_p=p^{rs}$, and since $|\A/M|_p=p^2$, we obtain that $sr=1$ or $2$.

We finish the proof by considering the two subcases: $H/M$ is insoluble or soluble.

\vskip 0.2cm

\noindent{\bf Subcase 2.1:} $H/M$ is insoluble.

In this subcase, $L_i/M$ are isomorphic non-abelain simple groups. We prove this subcase cannot happen by deriving contradictions. Recall that $sr=1$ or $2$.

Let $sr=1$. Then $|H/M|_p=p$, and therefore $|H|_p=p^2$. Since $H\unlhd \A$, $H\cap R(G)$ is a Sylow $p$-subgroup of $H$, implying $|H\cap R(G)|=p^2$, and then $R(G)\unlhd N$ yields that $H\cap R(G)$ is a normal subgroup of order $p^2$ in $N$, contrary to Lemma~\ref{arc-trans}.

Let $rs=2$. Then $|H/M|_p=p^2$ and $|H|_p=p^3$. This yields $R(G)\leq H$ and $H$ is transitive on $\Ga$, so $H/M$ is transitive on $V(\Ga_M)$. Note that $(r,s)=(1,2)$ or $(2,1)$.

For $(r,s)=(1,2)$, $H/M$ is a nonabelian simple group. By Propostion~\ref{Praeger}, $(H/M)_u$ for $u\in V(\Ga_M)$ is a $p'$-group because $H_1\leq \A_1$ is a $p'$-group, and by Proposition~\ref{2-trans}, $H/M$ is $2$-transitive on $V(\Ga_M)$, forcing that $\Ga_M$ is a complete group of order $p^2$, a contradiction.

For $(r,s)=(2,1)$, $H/M\cong L_1/M\times L_2/M$, where $L_1/M$ and $L_2/M$ are isomorphic nonabelain simple groups and $|L_i/M|_p=p$. It follows that $|H|_p=p^3$ and $|L_i|_p=p^2$ for $1\leq i\leq 2$. Since $H\unlhd \A$, we derive $R(G)\leq H$. Note that $H$ has no subgroup of index $2$. Since $\Sig$ is bipartite, it is $H$-semisymmetric. Let $\Delta_1$ and $\Delta_2$ be the two parts of $\Sig$. Then $|\Delta_1|=|\Delta_2|=p^2$, and $H$ is transitive on both $\Delta_1$ and $\Delta_2$.

Suppose $(L_1)_u=1$ for some $u\in V(\Sig)=\Delta_1\cup\Delta_2$. By Proposition~\ref{orbit-stab}, $|L_1|=|u^{L_1}|$, and since $L_1\unlhd H$ and $|\Delta_1|=|\Delta_2|=p^2$, we derive $|L_1|=p$ or $p^2$, contrary to the insolubleness of $L_1$. Thus $(L_1)_u\not=1$. Since $\Sig$ has prime valency $p$, $H_u$ is primitive on the neighbourhood $\Sig(u)$ of $u$ in $\Sig$, and since $(L_1)_u\unlhd H_u$, $(L_1)_u$ is transitive on $\Sig(u)$, which implies that $|(L_1)_u|_p=p$. Since $|L_1|_p=p^2$, each orbit of $L_1$ on $\Delta_1$ or $\Delta_2$ has length $p$.

Let $x\in \Delta_1$ and $y\in \Delta_2$ be adjacent in $\Sig$, and let $\Delta_{11}$ and $\Delta_{21}$ be the orbits of $L_1$ containing $x$ and $y$, respectively. Then $|\Delta_{11}|=|\Delta_{21}|=p$. Since $(L_1)_x$ is transitive on $\Sig(x)$, $x$ is adjacent to each vertex in $\Delta_{21}$, and therefore, each vertex in $\Delta_{11}$ is adjacent to each vertex in $\Delta_{21}$, that is, the induced subgroup $[\Delta_{11}\cup \Delta_{21}]$ is the complete bipartite graph $\K_{p,p}$. It follows that $\Sig\cong p\K_{p,p}$, contrary to the connectedness of $\Sig$.

\vskip 0.2cm
\noindent{\bf Subcase 2.2:} $H/M$ is soluble.

In this case, $|H|=p^2$ or $p^3$. Recall that $H\unlhd \A$. If $|H|=p^2$ then $H\leq R(G)$ and $N$ has normal subgroup of order $p^2$, contradicts Lemma~\ref{arc-trans}. Thus, $|H|=p^3$ and $\A$ has a normal Sylow $p$-subgroup, as required. This completes the proof. \qed

Now we are ready to finish the proof.

\vskip 0.2cm

{\noindent\bf Proof of Theorem~\ref{Thm-1}.}
By Lemmas~\ref{arc-trans} and \ref{normal}, $\Ga$ is a arc-transitive normal Cayley graph. In particular, $\Ga$ is $1$-distance transitive. Since $S=\{a^i,b^i\ |\ 1\leq i\leq p-1\}$, $\Ga$ has girth $3$, so it is not $2$-arc-transitive.

Recall that $G=\la a,b,c\mid a^p=b^p=c^p=1,[a,b]=c,[c,a]=[c,b]=1\ra$. Clearly,
\begin{eqnarray*}
&~&\Ga_1(1)=S=\{a^i,b^i \ |\ 1\leq i\leq p-1\}, \\
&~& \Ga_2(1)=\{b^ja^i,a^jb^i \ |\ 1\leq i,j\leq p-1\}.
\end{eqnarray*}

Note that $\Aut(G,S)=\la\a,\b,\g\mid \a^{p-1}=\b^{p-1}=\g^2=1, \a^\b=\a,\a^\g=\b\ra$, where $a^\a=a^t,b^\a=b,c^\a=c^t$, $a^\b=a,b^\b=b^t,c^\b=c^t$, $a^\g=b,b^\g=a$ and $c^\g=c^{-1}$. Then $(ba)^{\a^i\b^j}=b^{t^i}a^{t^j}$, and since $\ZZ_p^*=\langle t\rangle$, we obtain that $\langle \a,\b\rangle$ is transitive on the set $\{ b^ja^i\ | \ 1\leq i,j\leq p-1\}$. Similarly, $\langle \a,\b\rangle$ is transitive on $\{a^jb^i\ | \ 1\leq i,j\leq p-1\}$. Furthermore, $\g$ interchanges the two sets $\{ b^ja^i\ | \ 1\leq i,j\leq p-1\}$ and  $\{a^jb^i\ | \ 1\leq i,j\leq p-1\}$. It follows that $\Aut(G,S)$ is transitive on $\Ga_2(1)$ and hence $\Ga$ is $2$-distance transitive.

Noting that $ab=bac$, we have that $b^{-1}ab=ac\in \Ga_3(1)$ and  $aba=ba^2c\in \Ga_3(1)$. Also it is easy to see that $(ac)^{\Aut(G,S)}=(ac)^{\langle \a,\b,\g\rangle}=\{a^ic^j,b^ic^j\ |\ 1\leq i,j\leq p-1\}$. Now it is easy to see that $ba^2c\not\in (ac)^{\Aut(G,S)}$, and since $\A_1=\Aut(G,S)$ by Proposition~\ref{Aut-Cay}, $\Ga$ is not distance-transitive.   \qed

\noindent {\bf Proof of Corollary~\ref{cor}:} Recall that $\Sig$ is the clique graph of $\Ga$. By the first paragraph in the proof of Lemma~\ref{clique} and the definition of $\Cos(G, \la a\ra, \la b\ra)$ in Corollary~\ref{cor}, we have $\Sig=\Cos(G, \la a\ra, \la b\ra)$. Again by Lemma~\ref{clique},
$\Sig$ is $R(G)$-semisymmetric, and since $|E(\Sig)|=(2p^2\cdot p)/2=p^3=|R(G)|$, $R(G)$ is regular on the edge set $E(\Sig)$ of $\Sig$. Thus, the line graph of $\Sig$ is a Cayley graph on $G$.

For a given edge $\{\langle a\rangle x,\langle b\rangle y\}\in E(\Sig)$, we have $|\langle a\rangle x\cap \langle b\rangle y|=1$, and then we may identify this edge with the unique element in $\langle a\rangle x\cap \langle b\rangle y$. Note that $\Sig$ has valency $2(p-1)$. Then the edge $1=\langle a\rangle \cap \langle b\rangle$  in $\Sig$ is exactly incident to all edges in $S=\{a^i,b^i\ |\ 1\leq i\leq p-1\}$, because $\{a^i\}=\langle a\rangle\cap \langle b\rangle a^i$ and $\{b^i\}=\langle b\rangle\cap \langle a\rangle b^i$. It follows that $\Ga=\Cay(G,S)$ is exactly the line graph of $\Sig$.

If $\a\in\Aut(\Sig)$ fixes each edge in $\Sig$ then $\a$ fixes all vertices of $\Sig$, that is, $\Aut(\Sig)$ acts faithfully on $\Ga$. Thus, we may view $\Aut(\Sig)$ as a subgroup of $\Aut(\Ga)$. By Lemmas~\ref{clique} and \ref{normal}, we have $\Aut(\Ga)=\Aut(\Sig)=R(G)\rtimes\Aut(G,S)$.

Recall that $\Aut(G,S)=\langle \a,\b,\g\rangle$ and $\Sigma$ is arc-transitive. Since $a^\b=a$, $b^\b=b^t$ and $c^\b=c^t$, where $\ZZ_p^*=\langle t\rangle$, $\langle\b\rangle$ fixes the arc $(\langle a\rangle,\langle b\rangle)$ in $\Sigma$ and is transitive on the vertex set $\{\langle a\rangle b^i\ |\ 1\leq i\leq p-1 \}$, where $\{\langle a\rangle \}\cup \{\langle a\rangle b^i\ |\ 1\leq i\leq p-1 \}$ is the neighbourhood of $\langle b\rangle$ in $\Sig$. Thus, $\Sig$ is $2$-arc-transitive. Since $a^\a=a^t$, $b^\a=b$ and $c^\a=c^t$, $\langle\a\rangle$ fixes the $2$-arc $(\langle a\rangle,\langle b\rangle,\langle a\rangle b)$ and is transitive on  the vertex set $\{\langle b \rangle a^ib\ |\ 1\leq i\leq p-1 \}$, where $\{\langle b\rangle \}\cup \{\langle b\rangle a^ib\ |\ 1\leq i\leq p-1 \}$ is the neighbourhood of $\langle a\rangle b$ in $\Sig$. It follows that $\Sig$ is $3$-arc-transitive. It is easy to see that the number of $3$-arcs in $\Sig$ equals to $|A|=2p^3(p-1)^2$, $A$ is regular on the set of $3$-arcs of $\Sig$. \qed

\end{document}